\documentclass[12pt,preprint]{aastex}

\usepackage{amsmath,amssymb,graphics,graphicx}
\usepackage{aas_macros}
\usepackage{natbib}
\newcommand{\mathsym}[1]{{}}
\newcommand{\unicode}[1]{{}}

\shorttitle{Non-linear stability\dots smaller primary}
\shortauthors{B. Ishwar and J. P. Sharma}

\begin{document}


\title{Non-linear stability in photogravitational non-planar restricted three body problem with oblate smaller primary}


\author{B. Ishwar and J. P. Sharma}
\affil{University Department of Mathematics,B.R.A. Bihar university, Muzaffarpur-842001}
\email{ishwar\_bhola@hotmail.com}
%
%

%


\begin{abstract}
We have discussed non-linear stability in photogravitational non-planar restricted three body problem with oblate smaller primary. By photogravitational we mean that both primaries are radiating. We normalised the Hamiltonian using Lie transform as in Coppola and Rand (1989). We transformed the system into Birkhoff's normal form. Lie transforms reduce the system to an equivalent simpler system which is immediately solvable. Applying Arnold’s theorem, we have found non-linear stability criteria. We conclude that $L_6$ is stable. We plotted graphs for $(\omega_1 ,D_2).$  They are rectangular hyperbola. 
\end{abstract}


\keywords{Non-linear stability:Photogravitational:Non-planar:Oblate primary:RTBP}

\section{Introduction}
G. Hori(1966,1967) applied a theorem by Lie in canonical transformation to the theory of general perturbations. Theorem is applicable to such cases where the undisturbed portion of Hamiltonian depends on angular variable as well as momentum variables. A. Deprit(1969) introduced the concept of Lie series to the cases where the generating function itself depends explicitly on the small parameter. Lie transforms define naturally a class of canonical mappings in the form of power series in the small parameter. They reviewed how a Lie series defines a canonical mapping as a formal power series of a small parameter $\epsilon$, provided the generating function itself does not depend upon $\epsilon$. This restriction is overcome by introducing Lie transform. They showed that how they naturally define the canonical transformations contemplated by Von Zeipel's method. Canonical mappings defined by Lie transforms as formal power series of a small parameter constitute the natural ingredient of Transformation Theory applied to Hamiltonian systems. 
                   
Orbital stability of quasi-periodic motions in multidimensional Hamiltonian systems was studied by Sokolskii(1978). With some applications to the Birkhoff's normal form along with its generalized form by K. R. Meyer, the restricted problem of three bodies near $L_4$, the Birkhoff's normalization procedure, and the singular perturbation, of Hamiltonian systems have been  discussed by Liu(1985).   K.R. Meyer and D.S.Schmidt(1986) established the full stability of Lagrange equilibrium point in the planar restricted three body problem even in the case when $\mu= \mu_c$ . Hamiltonian is normalised up to order six and then KAM theory is applied. This establishes the stability of the equilibrium in degenerate case.  A.P. Markeev(1966) and K.T. Alfriend(1970,1971) have shown that $L_4$ is unstable when the mass ratio is equal to $\mu_2$ or $\mu_3$.   The Lie transform method is an efficient perturbation scheme which  explicitly generates the functional form of the reduced Hamiltonian under an implicitly defined canonical periodic near identity - transformation.

Coppola \& Rand(1989) applied a method of Lie Transforms, a perturbation method for differential equations to a general class of Hamiltonian systems using computer algebra. They developed explicit formulas for transforming the system into Birkhoff normal form. They formed explicit nonlinear stability criteria solely in terms of H for systems where the linear stability is               inconclusive. They applied these results to the non-linear stability of $L_4$ in the Circular RTBP. At $L_4$ , Arnold's theorem(1961) must be used since a Lyapunov function cannot be found. They confirmed the previous computations of Deprit and   Deprit –Bartholome (1967), Meyer and Schmidt (1986). Algorithms of linear and nonlinear normalization of a Hamiltonian system near an equilibrium point were described by Maciejewski, and Gozdziewski(1991). 

A. Jorba(1997) described the effective computation of normal forms,
 centre manifolds and first integrals in Hamiltonian Mechanics. These kind of conclusions are very useful. They allow, for example, to give explicit estimates on the diffusion time or to compute invariant tori. Their approach presented here is based on using algebraic manipulation for the formal series but taking numerical   co-efficients for them.  

 J. Palacian \& P. Yanguas (2000) described the reduction of perturbed Hamiltonian systems. They used a technique based on Lie-transformations. By extending an integral of the unperturbed part to the whole transformed system up to a certain order of approximation, the number of degrees of freedom of such a system is reduced, under certain conditions.  The idea of reducing a perturbed system is valid not only for Hamiltonians, but also for any system of differential equations. Recently there has been a resurgence of this subject. The problem of building formal integrals for Hamiltonian systems has received a wide treatment in the last forty five years. The results have been applied in fields such as Molecular Physics or Astrodynamics. For example, in galactic models of three degrees of freedom, the search for third integral becomes very important to analyze the onset of Chaos. Their approach consists in generalising the concept of normal forms by selecting a function $G(x)$ and proposing thereafter a symplectic change of variables.

 The nonlinear stability of triangular equilibrium points  was studied by Kushvah et al (2007) in the generalised photogravitational restricted three body
problem with Poynting-Robertson drag. They have performed first and
second order normalization of the Hamiltonian of the problem and applied KAM theorem to examine the condition of non-linear stability. After computation they 
have found three critical mass ratios and concluded that triangular points are stable in the nonlinear sense except three critical mass ratios at which KAM theorem fails.

 Hence, we thought to examine the Non-linear stability of $L_6$, equilibrium point of Non-planar photogravitational restricted three body problem with oblate smaller primary. We examined the linear stability of above problem in Shankaran(2011) where $q_1$ is radiation pressure of bigger primary and $q_2$ that of smaller. We have found that $L_6$ is unstable due to positive real part in complex roots. Now, we proceed to normalize the Hamiltonian of the problem as in Coppola and Rand (1989) using Lie transform. We transformed the system into Birkhoff's normal form. Using Arnold's theorem, we have found non-linear stability criteria in terms of H. Lie transforms reduce the system to an equivalent simpler system which is immediately solvable. Finally we find that $D_2 \neq 0.$ Hence, according to Arnold theorem, we conclude that $L_6$  is stable. We plotted graphs for $(\omega_1 ,D_2).$ They are rectangular hyperbola.

\section{Computation of $D_2$}
We suppose that $q$ is the radiation co-efficient of bigger primary and $Q$ that of smaller primary in photogravitational non-planar restricted three body problem. We want to apply Arnold's theorem(1961).  We applied method of Lie transform  using computer algebra as given in research article of Coppola and Rand(1989). We suppose $L_n=\{,w_n\}$ and $S_n=\frac{1}{n} \sum_{m=0}^{n-1} L_{n-m} S_m$,(where $n=1,2,3\dots$)  be the operators. Then we transformed the near-identity transformation  $(x_m,z_m)$ to $(X_n,Z_m)$ variables.  The Kamiltonian is given by 

\begin{equation}
K_n=H_n+\frac{1}{2} \{H_n,w_n\}+\frac{1}{2} \sum_{m=0}^{n-1}[ L_{n-m} K_m+mS_{n-m}H_m] \quad n=2,3,4\dots\end{equation} 
We choose the generating function $w_n$ in such a way so as to best simplify Kamiltonian and each term to be canceled will be of the form $A X_1^jY_1^lX_2^rY_2^s$ where $A$ is a constant. We choose $w_n$ to be a sum of terms, one for each term to be canceled, of the form  $w_n=A X_1^jY_1^lX_2^rY_2^s$ where $B = \frac{i A_n}{\omega_1(l-j)+\omega_3(s-r)}$ with $n=j+l+r+s-2.$  Here $\omega_1$ and $\omega_3$  are  the basic frequencies which are rationally independent and we suppose that the frequency in vertical direction is constant i.e. $\omega_3= 1$. If $j=l$ and $s=r$ then $B$ will be infinite hence the terms $ (X_1Y_1)^j(X_2Y_2)^r$. The $D_2$  is given as (please see Coppola and Rand(1989)) 

\begin{equation}
D_2=-(K2200\omega_3^2 +K1111\omega_1\omega_3+K0022\omega_1^2)
\end{equation} 
We performed computation using Mathematica and found the following results:

\begin{eqnarray}
\text{K0022}&&=\text{H0022}+i\left(\frac{1 }{\omega _1\text{  }}\text{H0111}\times\text{H1011}+\frac{1 }{\omega _3\text{  }}(\text{H0012}\times\text{H0021}\right.\nonumber\\&&+\left.\text{H0003}\times\text{H0030})\right.\nonumber\\&&+
\left.\frac{ 1 }{\left( \omega _1-2\omega _3\right)}\text{H0120}\times\text{H1002}+\frac{1}{\left(\omega _1+2\omega _3\right)}\text{H0102}\times \text{H1020}\right);
\end{eqnarray}
\begin{eqnarray}
\text{K1111}&&=\text{H1111}+i\left(\frac{ 1 }{ \omega _1 }(\text{H1011}\times\text{H1200}\right.\nonumber\\&&+\left.\text{H0111}\times\text{H2100})+\frac{ 1}{\left( \omega _1-2\omega
_3\right)} \text{H0120}\times\text{H1002}+\right.\nonumber\\&&
\frac{ 1 }{2 \omega _3 }(\text{H1101}\times\text{H0021}+\text{H1110}*\text{H0012})+\frac{ 1}{\left( 2\omega _1-\omega _3\right)}\text{H0210}\times\text{H2001}+\nonumber\\&&
\left.\frac{1}{\left(\omega _1+2\omega _3\right)}\text{H0102}\times\text{H1020}+\frac{1}{\left(2\omega _1+\omega _3\right)}\text{H0201}\times\text{H2010}\right);
\end{eqnarray}
\begin{eqnarray}
\text{K2200}&&=\text{H2200}+i\left(\frac{ 1}{\omega _3\text{  }}\text{H1101}\times\text{H1110}+\frac{1 }{\omega _1\text{  }}(\text{H1200}\times\text{H2100}\right.\\&&
+\text{H0300}\times\text{H3000})+\left.\frac{ 1}{\left(2 \omega _1-\omega _3\right)}\text{H0210}\times\text{H2001}+\nonumber\frac{ 1 }{\left(2 \omega _1+\omega _3\right)}\text{H0201}\times\text{H2010}\right)\nonumber\end{eqnarray}
\begin{eqnarray}
\text{K2200}&&=\frac{\left(5 a_1^2-6 b_1 \omega _1\right) \omega _3 \left(4 \omega _1^2-\omega _3^2\right)+2 a_2^2 \omega _1 \left(4 \omega _1^2+\omega_1 \omega _3-\omega _3^2\right)}{16 \omega _1^3 \omega _3-4 \omega _1 \omega _3^3}\\
\text{K1111}&&=\frac{1}{8} \left(-8 b_3+\frac{12 a_1 a_3}{\omega _1}+\frac{a_3^2}{\omega _1-2 \omega _3}+\frac{a_2^2}{2 \omega _1-\omega _3}+\frac{6a_2 a_4}{\omega _3}+\frac{a_2^2}{2 \omega _1+\omega _3}+\frac{a_3^2}{\omega _1+2 \omega _3}\right)\nonumber\\&&
\end{eqnarray}

\begin{eqnarray}
\text{K0022}&&= \frac{1}{8} \left(-12 b_5+\frac{10 a_4^2}{\omega _3}+a_3^2 \left(\frac{4}{\omega _1}+\frac{2 \omega _1}{\omega _1^2-4 \omega _3^2}\right)\right)\end{eqnarray}

\begin{equation}D_2=-\left(\text{K2200} \omega _3{}^2+ \text{K1111} \omega _1 \omega _3 +\text{K0022} \omega _1{}^2\right)\end{equation}

\noindent\(=\frac{1}{4} \left(-3 a_2 a_4 \omega _1+6 b_5 \omega _1^2-\frac{5 a_4^2 \omega _1^2}{\omega _3}-6 a_1 a_3 \omega _3+4 b_3 \omega
_1 \omega _3+\right.\\
\left.6 b_1 \omega _3^2-\frac{5 a_1^2 \omega _3^2}{\omega _1}-\frac{a_3^2 \omega _1 \left(3 \omega _1^2+\omega _1 \omega _3-8 \omega _3^2\right)}{\omega
_1^2-4 \omega _3^2}+\frac{2 a_2^2 \omega _3 \left(5 \omega _1^2+\omega _1 \omega _3-\omega _3^2\right)}{-4 \omega _1^2+\omega _3^2}\right)\)

where \\

\noindent\(a\text{=}(1-\mu )+\frac{6(3)^{1/2} ((1-\mu ))(1-q)A^{3/2}}{\mu  Q};\)

\noindent\(c\text{=}(3*A)^{1/2} -\frac{9 ((1-\mu ))q A^2}{\mu  Q};\)

\noindent\(a_1=\frac{-q \mu ^4+Q \sqrt{(-1+\mu )^2} \sqrt{\mu ^2}-2 Q \sqrt{(-1+\mu )^2} \mu  \sqrt{\mu ^2}+Q \sqrt{(-1+\mu )^2} \left(\mu ^2\right)^{3/2}}{(-1+\mu
)^2 \sqrt{(-1+\mu )^2} \mu ^4}-\\
\left(5 \left(-3 q \mu ^6+2 Q \sqrt{(-1+\mu )^2} \sqrt{\mu ^2}-8 Q \sqrt{(-1+\mu )^2} \mu  \sqrt{\mu ^2}-8 Q \sqrt{(-1+\mu )^2} \mu ^3 \sqrt{\mu
^2}+\right.\right.\\
\left.\left.\left.2 Q \sqrt{(-1+\mu )^2} \mu ^4 \sqrt{\mu ^2}+12 Q \sqrt{(-1+\mu )^2} \left(\mu ^2\right)^{3/2}\right) A\right)\right/\\
\left((-1+\mu )^4 \sqrt{(-1+\mu )^2} \mu ^6\right)+\left(24 \left(\sqrt{3} q \mu ^5-\sqrt{3} q^2 \mu ^5+\sqrt{3} Q \sqrt{(-1+\mu )^2} \sqrt{\mu
^2}-\right.\right.\\
\sqrt{3} q Q \sqrt{(-1+\mu )^2} \sqrt{\mu ^2}-3 \sqrt{3} Q \sqrt{(-1+\mu )^2} \mu  \sqrt{\mu ^2}+\\
3 \sqrt{3} q Q \sqrt{(-1+\mu )^2} \mu  \sqrt{\mu ^2}-\sqrt{3} Q \sqrt{(-1+\mu )^2} \mu ^3 \sqrt{\mu ^2}+\\
\sqrt{3} q Q \sqrt{(-1+\mu )^2} \mu ^3 \sqrt{\mu ^2}+3 \sqrt{3} Q \sqrt{(-1+\mu )^2} \left(\mu ^2\right)^{3/2}-\\
\left.\left.3 \sqrt{3} q Q \sqrt{(-1+\mu )^2} \left(\mu ^2\right)^{3/2}\right) A^{3/2}\right)/\left(Q (-1+\mu )^2 \sqrt{(-1+\mu )^2} \mu ^6\right)-\\
\left(945 \left(q \mu ^8+Q \sqrt{(-1+\mu )^2} \sqrt{\mu ^2}-6 Q \sqrt{(-1+\mu )^2} \mu  \sqrt{\mu ^2}-20 Q \sqrt{(-1+\mu )^2} \mu ^3 \sqrt{\mu
^2}+\right.\right.\\
15 Q \sqrt{(-1+\mu )^2} \mu ^4 \sqrt{\mu ^2}-6 Q \sqrt{(-1+\mu )^2} \mu ^5 \sqrt{\mu ^2}+\\
\left.\left.\left.Q \sqrt{(-1+\mu )^2} \mu ^6 \sqrt{\mu ^2}+15 Q \sqrt{(-1+\mu )^2} \left(\mu ^2\right)^{3/2}\right) A^2\right)\right/\\
\left(8 \left((-1+\mu )^6 \sqrt{(-1+\mu )^2} \mu ^8\right)\right)+O[A]^{5/2};\)

\noindent\( {a_2=\left(-\frac{6 \sqrt{3} q}{\left((1-\mu )^2\right)^{5/2}}+\frac{15 \sqrt{3} q \mu }{2 \left((1-\mu )^2\right)^{5/2}}-\frac{3
\sqrt{3} q \sqrt{(1-\mu )^2} \mu }{2 (1-\mu )^6}-\frac{6 \sqrt{3} Q \sqrt{\mu ^2}}{\mu ^5}\right) \sqrt{A}-}\\
 {\frac{135 \left(\sqrt{3} q\right) A^{3/2}}{2 \left((-1+\mu )^5 \sqrt{(-1+\mu )^2}\right)}+}\\
 {\left(54 \left(-10 q \mu ^6+9 q^2 \mu ^6+q^2 \mu ^7+10 Q \sqrt{(-1+\mu )^2} \sqrt{\mu ^2}-10 q Q \sqrt{(-1+\mu )^2} \sqrt{\mu ^2}-\right.\right.}\\
 {40 Q \sqrt{(-1+\mu )^2} \mu  \sqrt{\mu ^2}+39 q Q \sqrt{(-1+\mu )^2} \mu  \sqrt{\mu ^2}-}\\
 {40 Q \sqrt{(-1+\mu )^2} \mu ^3 \sqrt{\mu ^2}+34 q Q \sqrt{(-1+\mu )^2} \mu ^3 \sqrt{\mu ^2}+10 Q \sqrt{(-1+\mu )^2} }\\
 {\mu ^4 \sqrt{\mu ^2}-6 q Q \sqrt{(-1+\mu )^2} \mu ^4 \sqrt{\mu ^2}-q Q \sqrt{(-1+\mu )^2} \mu ^5 \sqrt{\mu ^2}+}\\
 {\left.\left.\left.60 Q \sqrt{(-1+\mu )^2} \left(\mu ^2\right)^{3/2}-56 q Q \sqrt{(-1+\mu )^2} \left(\mu ^2\right)^{3/2}\right) A^2\right)\right/}\\
 {\left(Q (-1+\mu )^3 \sqrt{(-1+\mu )^2} \mu ^7\right)+O[A]^{5/2};}\)

\noindent\( {a_3=\left(\frac{3 q (1-\mu )}{2 \left((1-\mu )^2\right)^{5/2}}-\frac{3 q (1-\mu ) \mu }{2 \left((1-\mu )^2\right)^{5/2}}-\frac{3
Q \sqrt{\mu ^2}}{2 \mu ^4}\right)+}\\
 {\left(-\frac{45 q (1-\mu )}{2 \left((1-\mu )^2\right)^{7/2}}-\frac{45 q}{4 (1-\mu ) \left((1-\mu )^2\right)^{5/2}}+\frac{45 q (1-\mu ) \mu }{2
\left((1-\mu )^2\right)^{7/2}}+\frac{45 q \mu }{4 (1-\mu ) \left((1-\mu )^2\right)^{5/2}}+\right.}\\
 {\left.\frac{45 Q \sqrt{\mu ^2}}{2 \mu ^6}\right) A+\left(\frac{36 \sqrt{3} (1-q) q (1-\mu )}{Q \left((1-\mu )^2\right)^{5/2}}-\frac{36 \sqrt{3}
(1-q) q (1-\mu )}{Q \left((1-\mu )^2\right)^{5/2} \mu }-\right.}\\
 {\left.\frac{36 \sqrt{3} \sqrt{\mu ^2}}{\mu ^6}+\frac{36 \sqrt{3} q \sqrt{\mu ^2}}{\mu ^6}+\frac{36 \sqrt{3} \sqrt{\mu ^2}}{\mu ^5}-\frac{36
\sqrt{3} q \sqrt{\mu ^2}}{\mu ^5}\right) A^{3/2}+}\\
 {\left(\frac{945 q}{4 (1-\mu ) \left((1-\mu )^2\right)^{7/2}}+\frac{945 q}{16 (1-\mu )^3 \left((1-\mu )^2\right)^{5/2}}-\frac{945 q \mu }{4 (1-\mu
) \left((1-\mu )^2\right)^{7/2}}-\right.}\\
 {\left.\frac{945 q \mu }{16 (1-\mu )^3 \left((1-\mu )^2\right)^{5/2}}+\frac{4725 Q \sqrt{\mu ^2}}{16 \mu ^8}\right) A^2+O[A]^{5/2};}\)

\noindent\( {a_4=\left(\frac{3 \sqrt{3} q}{2 \left((1-\mu )^2\right)^{5/2}}-\frac{3 \sqrt{3} q \sqrt{(1-\mu )^2} \mu }{2 (1-\mu )^6}+\frac{3 \sqrt{3}
Q \sqrt{\mu ^2}}{2 \mu ^5}\right) \sqrt{A}+}\\
 {\frac{75 \sqrt{3} q A^{3/2}}{4 (-1+\mu )^5 \sqrt{(-1+\mu )^2}}+\left(\frac{3}{2} q \left(\frac{\frac{9 q}{Q}-\frac{9 q}{Q \mu }}{\left((1-\mu
)^2\right)^{5/2}}-\frac{90 (1-q)}{Q \left((1-\mu )^2\right)^{5/2} \mu }\right)+\right.}\\
 {\frac{135 q \sqrt{(1-\mu )^2}}{Q (1-\mu )^6}-\frac{243 q^2 \sqrt{(1-\mu )^2}}{2 Q (1-\mu )^6}-\frac{27 q^2 \sqrt{(1-\mu )^2} \mu }{2 Q (1-\mu
)^6}+\frac{135 \sqrt{\mu ^2}}{\mu ^7}-}\\
 {\left.\frac{135 q \sqrt{\mu ^2}}{\mu ^7}-\frac{135 \sqrt{\mu ^2}}{\mu ^6}+\frac{243 q \sqrt{\mu ^2}}{2 \mu ^6}+\frac{27 q \sqrt{\mu ^2}}{2 \mu
^5}\right) A^2+O[A]^{5/2};}\)

\noindent\( {b_1=\left(-q \mu ^5-Q \sqrt{(-1+\mu )^2} \sqrt{\mu ^2}+3 Q \sqrt{(-1+\mu )^2} \mu  \sqrt{\mu ^2}+\right.}\\
 {\left.Q \sqrt{(-1+\mu )^2} \mu ^3 \sqrt{\mu ^2}-3 Q \sqrt{(-1+\mu )^2} \left(\mu ^2\right)^{3/2}\right)/\left((-1+\mu )^3 \sqrt{(-1+\mu )^2}
\mu ^5\right)-}\\
 {\left(15 \left(-3 q \mu ^7-2 Q \sqrt{(-1+\mu )^2} \sqrt{\mu ^2}+10 Q \sqrt{(-1+\mu )^2} \mu  \sqrt{\mu ^2}+20 Q \sqrt{(-1+\mu )^2} \right.\right.}\\
 {\mu ^3 \sqrt{\mu ^2}-10 Q \sqrt{(-1+\mu )^2} \mu ^4 \sqrt{\mu ^2}+2 Q \sqrt{(-1+\mu )^2} \mu ^5 \sqrt{\mu ^2}-}\\
 {\left.\left.20 Q \sqrt{(-1+\mu )^2} \left(\mu ^2\right)^{3/2}\right) A\right)/\left(2 \left((-1+\mu )^5 \sqrt{(-1+\mu )^2} \mu ^7\right)\right)+}\\
 {\left(30 \left(\sqrt{3} q \mu ^6-\sqrt{3} q^2 \mu ^6-\sqrt{3} Q \sqrt{(-1+\mu )^2} \sqrt{\mu ^2}+\sqrt{3} q Q \sqrt{(-1+\mu )^2} \sqrt{\mu ^2}+\right.\right.}\\
 {4 \sqrt{3} Q \sqrt{(-1+\mu )^2} \mu  \sqrt{\mu ^2}-4 \sqrt{3} q Q \sqrt{(-1+\mu )^2} \mu  \sqrt{\mu ^2}+}\\
 {4 \sqrt{3} Q \sqrt{(-1+\mu )^2} \mu ^3 \sqrt{\mu ^2}-4 \sqrt{3} q Q \sqrt{(-1+\mu )^2} \mu ^3 \sqrt{\mu ^2}-}\\
 {\sqrt{3} Q \sqrt{(-1+\mu )^2} \mu ^4 \sqrt{\mu ^2}+\sqrt{3} q Q \sqrt{(-1+\mu )^2} \mu ^4 \sqrt{\mu ^2}-}\\
 {\left.6 \sqrt{3} Q \sqrt{(-1+\mu )^2} \left(\mu ^2\right)^{3/2}+6 \sqrt{3} q Q \sqrt{(-1+\mu )^2} \left(\mu ^2\right)^{3/2}\right) }\\
 {\left.A^{3/2}\right)/\left(Q (-1+\mu )^3 \sqrt{(-1+\mu )^2} \mu ^7\right)-}\\
 {\left(945 \left(q \mu ^9-Q \sqrt{(-1+\mu )^2} \sqrt{\mu ^2}+7 Q \sqrt{(-1+\mu )^2} \mu  \sqrt{\mu ^2}+35 Q \sqrt{(-1+\mu )^2} \mu ^3 \sqrt{\mu
^2}-\right.\right.}\\
 {35 Q \sqrt{(-1+\mu )^2} \mu ^4 \sqrt{\mu ^2}+21 Q \sqrt{(-1+\mu )^2} \mu ^5 \sqrt{\mu ^2}-}\\
 {\left.7 Q \sqrt{(-1+\mu )^2} \mu ^6 \sqrt{\mu ^2}+Q \sqrt{(-1+\mu )^2} \mu ^7 \sqrt{\mu ^2}-21 Q \sqrt{(-1+\mu )^2} \left(\mu ^2\right)^{3/2}\right)
}\\
 {\left.A^2\right)/\left(4 \left((-1+\mu )^7 \sqrt{(-1+\mu )^2} \mu ^9\right)\right)+O[A]^{5/2};}\)

\noindent\( {b_3=}\\
 {\left(-\frac{3 q}{\left((1-\mu )^2\right)^{5/2}}+\frac{3 q \mu }{\left((1-\mu )^2\right)^{5/2}}+\frac{15 Q \sqrt{\mu ^2}}{4 \mu ^7}-\frac{15
Q (1-\mu )^2 \sqrt{\mu ^2}}{4 \mu ^7}-\frac{15 Q \sqrt{\mu ^2}}{2 \mu ^6}+\frac{3 Q \sqrt{\mu ^2}}{4 \mu ^5}\right)+}\\
 {\left(\frac{945 q}{8 \left((1-\mu )^2\right)^{7/2}}-\frac{45 q}{8 (1-\mu )^2 \left((1-\mu )^2\right)^{5/2}}-\frac{45 q \sqrt{(1-\mu )^2}}{4
(1-\mu )^8}-\frac{945 q \mu }{8 \left((1-\mu )^2\right)^{7/2}}+\right.}\\
 {\left.\frac{45 q \mu }{8 (1-\mu )^2 \left((1-\mu )^2\right)^{5/2}}+\frac{45 q \sqrt{(1-\mu )^2} \mu }{4 (1-\mu )^8}+\frac{135 Q \sqrt{\mu ^2}}{2
\mu ^7}\right) A-}\\
 {\left(90 \left(\sqrt{3} q \mu ^6-\sqrt{3} q^2 \mu ^6-\sqrt{3} Q \sqrt{(-1+\mu )^2} \sqrt{\mu ^2}+\sqrt{3} q Q \sqrt{(-1+\mu )^2} \sqrt{\mu ^2}+\right.\right.}\\
 {4 \sqrt{3} Q \sqrt{(-1+\mu )^2} \mu  \sqrt{\mu ^2}-4 \sqrt{3} q Q \sqrt{(-1+\mu )^2} \mu  \sqrt{\mu ^2}+4 \sqrt{3} Q }\\
 {\sqrt{(-1+\mu )^2} \mu ^3 \sqrt{\mu ^2}-4 \sqrt{3} q Q \sqrt{(-1+\mu )^2} \mu ^3 \sqrt{\mu ^2}-\sqrt{3} Q \sqrt{(-1+\mu )^2} }\\
 {\mu ^4 \sqrt{\mu ^2}+\sqrt{3} q Q \sqrt{(-1+\mu )^2} \mu ^4 \sqrt{\mu ^2}-6 \sqrt{3} Q \sqrt{(-1+\mu )^2} \left(\mu ^2\right)^{3/2}+}\\
 {\left.\left.6 \sqrt{3} q Q \sqrt{(-1+\mu )^2} \left(\mu ^2\right)^{3/2}\right) A^{3/2}\right)/\left(Q (-1+\mu )^3 \sqrt{(-1+\mu )^2} \mu ^7\right)+}\\
 {\left(4725 \left(q \mu ^9-Q \sqrt{(-1+\mu )^2} \sqrt{\mu ^2}+7 Q \sqrt{(-1+\mu )^2} \mu  \sqrt{\mu ^2}+35 Q \sqrt{(-1+\mu )^2} \mu ^3 \sqrt{\mu
^2}-\right.\right.}\\
 {35 Q \sqrt{(-1+\mu )^2} \mu ^4 \sqrt{\mu ^2}+21 Q \sqrt{(-1+\mu )^2} \mu ^5 \sqrt{\mu ^2}-}\\
 {\left.7 Q \sqrt{(-1+\mu )^2} \mu ^6 \sqrt{\mu ^2}+Q \sqrt{(-1+\mu )^2} \mu ^7 \sqrt{\mu ^2}-21 Q \sqrt{(-1+\mu )^2} \left(\mu ^2\right)^{3/2}\right)
}\\
 {\left.A^2\right)/\left(4 (-1+\mu )^7 \sqrt{(-1+\mu )^2} \mu ^9\right)+O[A]^{5/2};}\)

\noindent\( {b_5=}\\
 {\left(\frac{3 q}{8 \left((1-\mu )^2\right)^{5/2}}-\frac{3 q \mu }{8 \left((1-\mu )^2\right)^{5/2}}+\frac{3 Q \sqrt{\mu ^2}}{8 \mu ^5}\right)+\left(-\frac{45
q}{16 (1-\mu )^2 \left((1-\mu )^2\right)^{5/2}}-\frac{45 q \sqrt{(1-\mu )^2}}{4 (1-\mu )^8}+\right.}\\
 {\left.\frac{45 q \mu }{16 (1-\mu )^2 \left((1-\mu )^2\right)^{5/2}}+\frac{45 q \sqrt{(1-\mu )^2} \mu }{4 (1-\mu )^8}-\frac{75 Q \sqrt{\mu ^2}}{8
\mu ^7}\right) A+}\\
 {\left(\frac{45 \sqrt{3} (1-q) q}{4 Q \left((1-\mu )^2\right)^{5/2}}-\frac{45 \sqrt{3} (1-q) q}{4 Q \left((1-\mu )^2\right)^{5/2} \mu }+\frac{45
\sqrt{3} \sqrt{\mu ^2}}{4 \mu ^7}-\frac{45 \sqrt{3} q \sqrt{\mu ^2}}{4 \mu ^7}-\right.}\\
 {\left.\frac{45 \sqrt{3} \sqrt{\mu ^2}}{4 \mu ^6}+\frac{45 \sqrt{3} q \sqrt{\mu ^2}}{4 \mu ^6}\right) A^{3/2}+}\\
 {\left(\frac{945 q}{64 (1-\mu )^4 \left((1-\mu )^2\right)^{5/2}}+\frac{315 q \sqrt{(1-\mu )^2}}{2 (1-\mu )^{10}}-\frac{945 q \mu }{64 (1-\mu
)^4 \left((1-\mu )^2\right)^{5/2}}-\right.}\\
 {\left.\frac{315 q \sqrt{(1-\mu )^2} \mu }{2 (1-\mu )^{10}}-\frac{11025 Q \sqrt{\mu ^2}}{64 \mu ^9}\right) A^2+O[A]^{5/2};}\)

From above results  $D_2$ upto first order  in  $A$ is

\noindent\( {D_2=\frac{9 q \sqrt{(1-\mu )^2} \omega _1^2}{16 (1-\mu )^6}-\frac{9 q \sqrt{(1-\mu )^2} \mu  \omega _1^2}{16 (1-\mu )^6}+\frac{9 Q \sqrt{\mu
^2} \omega _1^2}{16 \mu ^5}+}\\
 {\frac{9 q^2 \sqrt{(1-\mu )^2} \sqrt{(-1+\mu )^2} \omega _3}{4 (1-\mu )^6 (-1+\mu )^4}+\frac{9 Q^2 \omega _3}{4 (-1+\mu )^2 \mu ^6}-\frac{9 Q^2
\omega _3}{2 (-1+\mu )^2 \mu ^5}+\frac{9 Q^2 \omega _3}{4 (-1+\mu )^2 \mu ^4}-}\\
 {\frac{9 q^2 \sqrt{(1-\mu )^2} \sqrt{(-1+\mu )^2} \mu  \omega _3}{2 (1-\mu )^6 (-1+\mu )^4}+\frac{9 q^2 \sqrt{(1-\mu )^2} \sqrt{(-1+\mu )^2}
\mu ^2 \omega _3}{4 (1-\mu )^6 (-1+\mu )^4}-}\\
 {\frac{9 q Q \sqrt{(1-\mu )^2} \sqrt{\mu ^2} \omega _3}{4 (1-\mu )^6 (-1+\mu )^2}-\frac{9 q Q \sqrt{(1-\mu )^2} \sqrt{\mu ^2} \omega _3}{4 (1-\mu
)^6 (-1+\mu )^2 \mu ^4}-\frac{9 q Q \sqrt{(-1+\mu )^2} \sqrt{\mu ^2} \omega _3}{4 (-1+\mu )^4 \mu ^4}+}\\
 {\frac{9 q Q \sqrt{(1-\mu )^2} \sqrt{\mu ^2} \omega _3}{(1-\mu )^6 (-1+\mu )^2 \mu ^3}-\frac{27 q Q \sqrt{(1-\mu )^2} \sqrt{\mu ^2} \omega _3}{2
(1-\mu )^6 (-1+\mu )^2 \mu ^2}+\frac{9 q Q \sqrt{(1-\mu )^2} \sqrt{\mu ^2} \omega _3}{(1-\mu )^6 (-1+\mu )^2 \mu }-}\\
 {\frac{3 q \sqrt{(1-\mu )^2} \omega _1 \omega _3}{(1-\mu )^6}+\frac{3 q \sqrt{(1-\mu )^2} \mu  \omega _1 \omega _3}{(1-\mu )^6}-\frac{3 Q \sqrt{\mu
^2} \omega _1 \omega _3}{\mu ^5}-}\\
 {\frac{3 q \sqrt{(-1+\mu )^2} \omega _3^2}{2 (-1+\mu )^5}-\frac{3 Q \sqrt{\mu ^2} \omega _3^2}{2 (-1+\mu )^3 \mu ^5}+\frac{9 Q \sqrt{\mu ^2}
\omega _3^2}{2 (-1+\mu )^3 \mu ^4}-\frac{9 Q \sqrt{\mu ^2} \omega _3^2}{2 (-1+\mu )^3 \mu ^3}+}\\
 {\frac{3 Q \sqrt{\mu ^2} \omega _3^2}{2 (-1+\mu )^3 \mu ^2}-\frac{5 q^2 \omega _3^2}{4 (-1+\mu )^6 \omega _1}-\frac{5 Q^2 \omega _3^2}{4 (-1+\mu
)^4 \mu ^6 \omega _1}+\frac{5 Q^2 \omega _3^2}{(-1+\mu )^4 \mu ^5 \omega _1}-}\\
 {\frac{15 Q^2 \omega _3^2}{2 (-1+\mu )^4 \mu ^4 \omega _1}+\frac{5 Q^2 \omega _3^2}{(-1+\mu )^4 \mu ^3 \omega _1}-\frac{5 Q^2 \omega _3^2}{4
(-1+\mu )^4 \mu ^2 \omega _1}+\frac{5 q Q \sqrt{(-1+\mu )^2} \sqrt{\mu ^2} \omega _3^2}{2 (-1+\mu )^6 \mu ^4 \omega _1}-}\\
 {\frac{5 q Q \sqrt{(-1+\mu )^2} \sqrt{\mu ^2} \omega _3^2}{(-1+\mu )^6 \mu ^3 \omega _1}+\frac{5 q Q \sqrt{(-1+\mu )^2} \sqrt{\mu ^2} \omega
_3^2}{2 (-1+\mu )^6 \mu ^2 \omega _1}-\frac{27 q^2 \omega _1^3}{16 (1-\mu )^{10} \left(\omega _1^2-4 \omega _3^2\right)}-}\\
 {\frac{27 Q^2 \omega _1^3}{16 \mu ^6 \left(\omega _1^2-4 \omega _3^2\right)}+\frac{27 q^2 \mu  \omega _1^3}{4 (1-\mu )^{10} \left(\omega _1^2-4
\omega _3^2\right)}-\frac{81 q^2 \mu ^2 \omega _1^3}{8 (1-\mu )^{10} \left(\omega _1^2-4 \omega _3^2\right)}+}\\
 {\frac{27 q^2 \mu ^3 \omega _1^3}{4 (1-\mu )^{10} \left(\omega _1^2-4 \omega _3^2\right)}-\frac{27 q^2 \mu ^4 \omega _1^3}{16 (1-\mu )^{10} \left(\omega
_1^2-4 \omega _3^2\right)}+\frac{27 q Q \sqrt{(1-\mu )^2} \sqrt{\mu ^2} \omega _1^3}{8 (1-\mu )^6 \mu ^4 \left(\omega _1^2-4 \omega _3^2\right)}-}\\
 {\frac{27 q Q \sqrt{(1-\mu )^2} \sqrt{\mu ^2} \omega _1^3}{4 (1-\mu )^6 \mu ^3 \left(\omega _1^2-4 \omega _3^2\right)}+\frac{27 q Q \sqrt{(1-\mu
)^2} \sqrt{\mu ^2} \omega _1^3}{8 (1-\mu )^6 \mu ^2 \left(\omega _1^2-4 \omega _3^2\right)}-\frac{9 q^2 \omega _1^2 \omega _3}{16 (1-\mu )^{10} \left(\omega
_1^2-4 \omega _3^2\right)}-}\\
 {\frac{9 Q^2 \omega _1^2 \omega _3}{16 \mu ^6 \left(\omega _1^2-4 \omega _3^2\right)}+\frac{9 q^2 \mu  \omega _1^2 \omega _3}{4 (1-\mu )^{10}
\left(\omega _1^2-4 \omega _3^2\right)}-\frac{27 q^2 \mu ^2 \omega _1^2 \omega _3}{8 (1-\mu )^{10} \left(\omega _1^2-4 \omega _3^2\right)}+}\\
 {\frac{9 q^2 \mu ^3 \omega _1^2 \omega _3}{4 (1-\mu )^{10} \left(\omega _1^2-4 \omega _3^2\right)}-\frac{9 q^2 \mu ^4 \omega _1^2 \omega _3}{16
(1-\mu )^{10} \left(\omega _1^2-4 \omega _3^2\right)}+\frac{9 q Q \sqrt{(1-\mu )^2} \sqrt{\mu ^2} \omega _1^2 \omega _3}{8 (1-\mu )^6 \mu ^4 \left(\omega
_1^2-4 \omega _3^2\right)}-}\\
 {\frac{9 q Q \sqrt{(1-\mu )^2} \sqrt{\mu ^2} \omega _1^2 \omega _3}{4 (1-\mu )^6 \mu ^3 \left(\omega _1^2-4 \omega _3^2\right)}+\frac{9 q Q \sqrt{(1-\mu
)^2} \sqrt{\mu ^2} \omega _1^2 \omega _3}{8 (1-\mu )^6 \mu ^2 \left(\omega _1^2-4 \omega _3^2\right)}+\frac{9 q^2 \omega _1 \omega _3^2}{2 (1-\mu
)^{10} \left(\omega _1^2-4 \omega _3^2\right)}+}\\
 {\frac{9 Q^2 \omega _1 \omega _3^2}{2 \mu ^6 \left(\omega _1^2-4 \omega _3^2\right)}-\frac{18 q^2 \mu  \omega _1 \omega _3^2}{(1-\mu )^{10} \left(\omega
_1^2-4 \omega _3^2\right)}+\frac{27 q^2 \mu ^2 \omega _1 \omega _3^2}{(1-\mu )^{10} \left(\omega _1^2-4 \omega _3^2\right)}-}\\
 {\frac{18 q^2 \mu ^3 \omega _1 \omega _3^2}{(1-\mu )^{10} \left(\omega _1^2-4 \omega _3^2\right)}+\frac{9 q^2 \mu ^4 \omega _1 \omega _3^2}{2
(1-\mu )^{10} \left(\omega _1^2-4 \omega _3^2\right)}-\frac{9 q Q \sqrt{(1-\mu )^2} \sqrt{\mu ^2} \omega _1 \omega _3^2}{(1-\mu )^6 \mu ^4 \left(\omega
_1^2-4 \omega _3^2\right)}+}\\
 {\frac{18 q Q \sqrt{(1-\mu )^2} \sqrt{\mu ^2} \omega _1 \omega _3^2}{(1-\mu )^6 \mu ^3 \left(\omega _1^2-4 \omega _3^2\right)}-\frac{9 q Q \sqrt{(1-\mu
)^2} \sqrt{\mu ^2} \omega _1 \omega _3^2}{(1-\mu )^6 \mu ^2 \left(\omega _1^2-4 \omega _3^2\right)}+}\\
 {A \left(\frac{81 q^2 \omega _1}{4 (1-\mu )^{10}}+\frac{81 Q^2 \omega _1}{4 \mu ^8}-\frac{81 q^2 \mu  \omega _1}{2 (1-\mu )^{10}}+\frac{81 q^2
\mu ^2 \omega _1}{4 (1-\mu )^{10}}+\frac{81 q Q \sqrt{(1-\mu )^2} \sqrt{\mu ^2} \omega _1}{2 (1-\mu )^6 \mu ^5}-\right.}\\
 {\frac{81 q Q \sqrt{(1-\mu )^2} \sqrt{\mu ^2} \omega _1}{2 (1-\mu )^6 \mu ^4}-\frac{675 q \sqrt{(1-\mu )^2} \omega _1^2}{32 (1-\mu )^8}+\frac{675
q \sqrt{(1-\mu )^2} \mu  \omega _1^2}{32 (1-\mu )^8}-}\\
 {\frac{225 Q \sqrt{\mu ^2} \omega _1^2}{16 \mu ^7}-\frac{135 q^2 \omega _1^2}{16 (1-\mu )^{10} \omega _3}-\frac{135 Q^2 \omega _1^2}{16 \mu ^8
\omega _3}+\frac{135 q^2 \mu  \omega _1^2}{8 (1-\mu )^{10} \omega _3}-}\\
 {\frac{135 q^2 \mu ^2 \omega _1^2}{16 (1-\mu )^{10} \omega _3}-\frac{135 q Q \sqrt{(1-\mu )^2} \sqrt{\mu ^2} \omega _1^2}{8 (1-\mu )^6 \mu ^5
\omega _3}+\frac{135 q Q \sqrt{(1-\mu )^2} \sqrt{\mu ^2} \omega _1^2}{8 (1-\mu )^6 \mu ^4 \omega _3}-}\\
 {\frac{135 q^2 \sqrt{(1-\mu )^2} \sqrt{(-1+\mu )^2} \omega _3}{4 (1-\mu )^6 (-1+\mu )^6}-\frac{135 q^2 \sqrt{(1-\mu )^2} \sqrt{(-1+\mu )^2} \omega
_3}{4 (1-\mu )^8 (-1+\mu )^4}-}\\
 {\frac{135 q^2 \sqrt{(1-\mu )^2} \sqrt{(-1+\mu )^2} \omega _3}{8 (1-\mu )^7 (-1+\mu )^4}-\frac{45 Q^2 \omega _3}{2 (-1+\mu )^4 \mu ^8}-\frac{135
Q^2 \omega _3}{4 (-1+\mu )^2 \mu ^8}+\frac{90 Q^2 \omega _3}{(-1+\mu )^4 \mu ^7}+}\\
 {\frac{135 Q^2 \omega _3}{2 (-1+\mu )^2 \mu ^7}-\frac{135 Q^2 \omega _3}{(-1+\mu )^4 \mu ^6}-\frac{135 Q^2 \omega _3}{4 (-1+\mu )^2 \mu ^6}+\frac{90
Q^2 \omega _3}{(-1+\mu )^4 \mu ^5}-\frac{45 Q^2 \omega _3}{2 (-1+\mu )^4 \mu ^4}+}\\
 {\frac{135 q^2 \sqrt{(1-\mu )^2} \sqrt{(-1+\mu )^2} \mu  \omega _3}{2 (1-\mu )^6 (-1+\mu )^6}+\frac{135 q^2 \sqrt{(1-\mu )^2} \sqrt{(-1+\mu )^2}
\mu  \omega _3}{2 (1-\mu )^8 (-1+\mu )^4}+}\\
 {\frac{135 q^2 \sqrt{(1-\mu )^2} \sqrt{(-1+\mu )^2} \mu  \omega _3}{8 (1-\mu )^7 (-1+\mu )^4}-\frac{135 q^2 \sqrt{(1-\mu )^2} \sqrt{(-1+\mu )^2}
\mu ^2 \omega _3}{4 (1-\mu )^6 (-1+\mu )^6}-}\\
 {\frac{135 q^2 \sqrt{(1-\mu )^2} \sqrt{(-1+\mu )^2} \mu ^2 \omega _3}{4 (1-\mu )^8 (-1+\mu )^4}+\frac{45 q Q \sqrt{(1-\mu )^2} \sqrt{\mu ^2}
\omega _3}{2 (1-\mu )^6 (-1+\mu )^4}+}\\
 {\frac{135 q Q \sqrt{(1-\mu )^2} \sqrt{\mu ^2} \omega _3}{4 (1-\mu )^8 (-1+\mu )^2}+\frac{45 q Q \sqrt{(1-\mu )^2} \sqrt{\mu ^2} \omega _3}{2
(1-\mu )^6 (-1+\mu )^4 \mu ^6}+}\\
 {\frac{135 q Q \sqrt{(-1+\mu )^2} \sqrt{\mu ^2} \omega _3}{4 (-1+\mu )^4 \mu ^6}-\frac{135 q Q \sqrt{(1-\mu )^2} \sqrt{\mu ^2} \omega _3}{(1-\mu
)^6 (-1+\mu )^4 \mu ^5}+}\\
 {\frac{675 q Q \sqrt{(1-\mu )^2} \sqrt{\mu ^2} \omega _3}{2 (1-\mu )^6 (-1+\mu )^4 \mu ^4}+\frac{135 q Q \sqrt{(1-\mu )^2} \sqrt{\mu ^2} \omega
_3}{4 (1-\mu )^8 (-1+\mu )^2 \mu ^4}+}\\
 {\frac{135 q Q \sqrt{(1-\mu )^2} \sqrt{\mu ^2} \omega _3}{8 (1-\mu )^7 (-1+\mu )^2 \mu ^4}+\frac{135 q Q \sqrt{(-1+\mu )^2} \sqrt{\mu ^2} \omega
_3}{4 (-1+\mu )^6 \mu ^4}-}\\
 {\frac{450 q Q \sqrt{(1-\mu )^2} \sqrt{\mu ^2} \omega _3}{(1-\mu )^6 (-1+\mu )^4 \mu ^3}-\frac{135 q Q \sqrt{(1-\mu )^2} \sqrt{\mu ^2} \omega
_3}{(1-\mu )^8 (-1+\mu )^2 \mu ^3}-}\\
 {\frac{405 q Q \sqrt{(1-\mu )^2} \sqrt{\mu ^2} \omega _3}{8 (1-\mu )^7 (-1+\mu )^2 \mu ^3}+\frac{675 q Q \sqrt{(1-\mu )^2} \sqrt{\mu ^2} \omega
_3}{2 (1-\mu )^6 (-1+\mu )^4 \mu ^2}+}\\
 {\frac{405 q Q \sqrt{(1-\mu )^2} \sqrt{\mu ^2} \omega _3}{2 (1-\mu )^8 (-1+\mu )^2 \mu ^2}+\frac{405 q Q \sqrt{(1-\mu )^2} \sqrt{\mu ^2} \omega
_3}{8 (1-\mu )^7 (-1+\mu )^2 \mu ^2}-}\\
 {\frac{135 q Q \sqrt{(1-\mu )^2} \sqrt{\mu ^2} \omega _3}{(1-\mu )^6 (-1+\mu )^4 \mu }-\frac{135 q Q \sqrt{(1-\mu )^2} \sqrt{\mu ^2} \omega _3}{(1-\mu
)^8 (-1+\mu )^2 \mu }-}\\
 {\frac{135 q Q \sqrt{(1-\mu )^2} \sqrt{\mu ^2} \omega _3}{8 (1-\mu )^7 (-1+\mu )^2 \mu }+\frac{405 q \sqrt{(1-\mu )^2} \omega _1 \omega _3}{4
(1-\mu )^8}-\frac{405 q \sqrt{(1-\mu )^2} \mu  \omega _1 \omega _3}{4 (1-\mu )^8}+}\\
 {\frac{135 Q \sqrt{\mu ^2} \omega _1 \omega _3}{2 \mu ^7}+\frac{135 q \sqrt{(-1+\mu )^2} \omega _3^2}{4 (-1+\mu )^7}+\frac{45 Q \sqrt{\mu ^2}
\omega _3^2}{2 (-1+\mu )^5 \mu ^7}-\frac{225 Q \sqrt{\mu ^2} \omega _3^2}{2 (-1+\mu )^5 \mu ^6}+}\\
 {\frac{225 Q \sqrt{\mu ^2} \omega _3^2}{(-1+\mu )^5 \mu ^5}-\frac{225 Q \sqrt{\mu ^2} \omega _3^2}{(-1+\mu )^5 \mu ^4}+\frac{225 Q \sqrt{\mu
^2} \omega _3^2}{2 (-1+\mu )^5 \mu ^3}-\frac{45 Q \sqrt{\mu ^2} \omega _3^2}{2 (-1+\mu )^5 \mu ^2}+}\\
 {\frac{75 q^2 \omega _3^2}{2 (-1+\mu )^8 \omega _1}+\frac{25 Q^2 \omega _3^2}{(-1+\mu )^6 \mu ^8 \omega _1}-\frac{150 Q^2 \omega _3^2}{(-1+\mu
)^6 \mu ^7 \omega _1}+\frac{375 Q^2 \omega _3^2}{(-1+\mu )^6 \mu ^6 \omega _1}-}\\
 {\frac{500 Q^2 \omega _3^2}{(-1+\mu )^6 \mu ^5 \omega _1}+\frac{375 Q^2 \omega _3^2}{(-1+\mu )^6 \mu ^4 \omega _1}-\frac{150 Q^2 \omega _3^2}{(-1+\mu
)^6 \mu ^3 \omega _1}+\frac{25 Q^2 \omega _3^2}{(-1+\mu )^6 \mu ^2 \omega _1}-}\\
 {\frac{25 q Q \sqrt{(-1+\mu )^2} \sqrt{\mu ^2} \omega _3^2}{(-1+\mu )^8 \mu ^6 \omega _1}+\frac{100 q Q \sqrt{(-1+\mu )^2} \sqrt{\mu ^2} \omega
_3^2}{(-1+\mu )^8 \mu ^5 \omega _1}-}\\
 {\frac{375 q Q \sqrt{(-1+\mu )^2} \sqrt{\mu ^2} \omega _3^2}{2 (-1+\mu )^8 \mu ^4 \omega _1}+\frac{175 q Q \sqrt{(-1+\mu )^2} \sqrt{\mu ^2} \omega
_3^2}{(-1+\mu )^8 \mu ^3 \omega _1}-}\\
 {\frac{125 q Q \sqrt{(-1+\mu )^2} \sqrt{\mu ^2} \omega _3^2}{2 (-1+\mu )^8 \mu ^2 \omega _1}+\frac{405 q^2 \omega _1^3}{8 (1-\mu )^{12} \left(\omega
_1^2-4 \omega _3^2\right)}+\frac{405 q^2 \omega _1^3}{16 (1-\mu )^{11} \left(\omega _1^2-4 \omega _3^2\right)}+}\\
 {\frac{405 Q^2 \omega _1^3}{8 \mu ^8 \left(\omega _1^2-4 \omega _3^2\right)}-\frac{405 q^2 \mu  \omega _1^3}{2 (1-\mu )^{12} \left(\omega _1^2-4
\omega _3^2\right)}-\frac{1215 q^2 \mu  \omega _1^3}{16 (1-\mu )^{11} \left(\omega _1^2-4 \omega _3^2\right)}+}\\
 {\frac{1215 q^2 \mu ^2 \omega _1^3}{4 (1-\mu )^{12} \left(\omega _1^2-4 \omega _3^2\right)}+\frac{1215 q^2 \mu ^2 \omega _1^3}{16 (1-\mu )^{11}
\left(\omega _1^2-4 \omega _3^2\right)}-\frac{405 q^2 \mu ^3 \omega _1^3}{2 (1-\mu )^{12} \left(\omega _1^2-4 \omega _3^2\right)}-}\\
 {\frac{405 q^2 \mu ^3 \omega _1^3}{16 (1-\mu )^{11} \left(\omega _1^2-4 \omega _3^2\right)}+\frac{405 q^2 \mu ^4 \omega _1^3}{8 (1-\mu )^{12}
\left(\omega _1^2-4 \omega _3^2\right)}-\frac{405 q Q \sqrt{(1-\mu )^2} \sqrt{\mu ^2} \omega _1^3}{8 (1-\mu )^6 \mu ^6 \left(\omega _1^2-4 \omega
_3^2\right)}+}\\
 {\frac{405 q Q \sqrt{(1-\mu )^2} \sqrt{\mu ^2} \omega _1^3}{4 (1-\mu )^6 \mu ^5 \left(\omega _1^2-4 \omega _3^2\right)}-\frac{405 q Q \sqrt{(1-\mu
)^2} \sqrt{\mu ^2} \omega _1^3}{8 (1-\mu )^8 \mu ^4 \left(\omega _1^2-4 \omega _3^2\right)}-}\\
 {\frac{405 q Q \sqrt{(1-\mu )^2} \sqrt{\mu ^2} \omega _1^3}{16 (1-\mu )^7 \mu ^4 \left(\omega _1^2-4 \omega _3^2\right)}-\frac{405 q Q \sqrt{(1-\mu
)^2} \sqrt{\mu ^2} \omega _1^3}{8 (1-\mu )^6 \mu ^4 \left(\omega _1^2-4 \omega _3^2\right)}+}\\
 {\frac{405 q Q \sqrt{(1-\mu )^2} \sqrt{\mu ^2} \omega _1^3}{4 (1-\mu )^8 \mu ^3 \left(\omega _1^2-4 \omega _3^2\right)}+\frac{405 q Q \sqrt{(1-\mu
)^2} \sqrt{\mu ^2} \omega _1^3}{16 (1-\mu )^7 \mu ^3 \left(\omega _1^2-4 \omega _3^2\right)}-}\\
 {\frac{405 q Q \sqrt{(1-\mu )^2} \sqrt{\mu ^2} \omega _1^3}{8 (1-\mu )^8 \mu ^2 \left(\omega _1^2-4 \omega _3^2\right)}+\frac{135 q^2 \omega
_1^2 \omega _3}{8 (1-\mu )^{12} \left(\omega _1^2-4 \omega _3^2\right)}+\frac{135 q^2 \omega _1^2 \omega _3}{16 (1-\mu )^{11} \left(\omega _1^2-4
\omega _3^2\right)}+}\\
 {\frac{135 Q^2 \omega _1^2 \omega _3}{8 \mu ^8 \left(\omega _1^2-4 \omega _3^2\right)}-\frac{135 q^2 \mu  \omega _1^2 \omega _3}{2 (1-\mu )^{12}
\left(\omega _1^2-4 \omega _3^2\right)}-\frac{405 q^2 \mu  \omega _1^2 \omega _3}{16 (1-\mu )^{11} \left(\omega _1^2-4 \omega _3^2\right)}+}\\
 {\frac{405 q^2 \mu ^2 \omega _1^2 \omega _3}{4 (1-\mu )^{12} \left(\omega _1^2-4 \omega _3^2\right)}+\frac{405 q^2 \mu ^2 \omega _1^2 \omega
_3}{16 (1-\mu )^{11} \left(\omega _1^2-4 \omega _3^2\right)}-\frac{135 q^2 \mu ^3 \omega _1^2 \omega _3}{2 (1-\mu )^{12} \left(\omega _1^2-4 \omega
_3^2\right)}-}\\
 {\frac{135 q^2 \mu ^3 \omega _1^2 \omega _3}{16 (1-\mu )^{11} \left(\omega _1^2-4 \omega _3^2\right)}+\frac{135 q^2 \mu ^4 \omega _1^2 \omega
_3}{8 (1-\mu )^{12} \left(\omega _1^2-4 \omega _3^2\right)}-\frac{135 q Q \sqrt{(1-\mu )^2} \sqrt{\mu ^2} \omega _1^2 \omega _3}{8 (1-\mu )^6 \mu
^6 \left(\omega _1^2-4 \omega _3^2\right)}+}\\
 {\frac{135 q Q \sqrt{(1-\mu )^2} \sqrt{\mu ^2} \omega _1^2 \omega _3}{4 (1-\mu )^6 \mu ^5 \left(\omega _1^2-4 \omega _3^2\right)}-\frac{135 q
Q \sqrt{(1-\mu )^2} \sqrt{\mu ^2} \omega _1^2 \omega _3}{8 (1-\mu )^8 \mu ^4 \left(\omega _1^2-4 \omega _3^2\right)}-}\\
 {\frac{135 q Q \sqrt{(1-\mu )^2} \sqrt{\mu ^2} \omega _1^2 \omega _3}{16 (1-\mu )^7 \mu ^4 \left(\omega _1^2-4 \omega _3^2\right)}-\frac{135
q Q \sqrt{(1-\mu )^2} \sqrt{\mu ^2} \omega _1^2 \omega _3}{8 (1-\mu )^6 \mu ^4 \left(\omega _1^2-4 \omega _3^2\right)}+}\\
 {\frac{135 q Q \sqrt{(1-\mu )^2} \sqrt{\mu ^2} \omega _1^2 \omega _3}{4 (1-\mu )^8 \mu ^3 \left(\omega _1^2-4 \omega _3^2\right)}+\frac{135 q
Q \sqrt{(1-\mu )^2} \sqrt{\mu ^2} \omega _1^2 \omega _3}{16 (1-\mu )^7 \mu ^3 \left(\omega _1^2-4 \omega _3^2\right)}-}\\
 {\frac{135 q Q \sqrt{(1-\mu )^2} \sqrt{\mu ^2} \omega _1^2 \omega _3}{8 (1-\mu )^8 \mu ^2 \left(\omega _1^2-4 \omega _3^2\right)}-\frac{135 q^2
\omega _1 \omega _3^2}{(1-\mu )^{12} \left(\omega _1^2-4 \omega _3^2\right)}-\frac{135 q^2 \omega _1 \omega _3^2}{2 (1-\mu )^{11} \left(\omega _1^2-4
\omega _3^2\right)}-}\\
 {\frac{135 Q^2 \omega _1 \omega _3^2}{\mu ^8 \left(\omega _1^2-4 \omega _3^2\right)}+\frac{540 q^2 \mu  \omega _1 \omega _3^2}{(1-\mu )^{12}
\left(\omega _1^2-4 \omega _3^2\right)}+\frac{405 q^2 \mu  \omega _1 \omega _3^2}{2 (1-\mu )^{11} \left(\omega _1^2-4 \omega _3^2\right)}-}\\
 {\frac{810 q^2 \mu ^2 \omega _1 \omega _3^2}{(1-\mu )^{12} \left(\omega _1^2-4 \omega _3^2\right)}-\frac{405 q^2 \mu ^2 \omega _1 \omega _3^2}{2
(1-\mu )^{11} \left(\omega _1^2-4 \omega _3^2\right)}+\frac{540 q^2 \mu ^3 \omega _1 \omega _3^2}{(1-\mu )^{12} \left(\omega _1^2-4 \omega _3^2\right)}+}\\
 {\frac{135 q^2 \mu ^3 \omega _1 \omega _3^2}{2 (1-\mu )^{11} \left(\omega _1^2-4 \omega _3^2\right)}-\frac{135 q^2 \mu ^4 \omega _1 \omega _3^2}{(1-\mu
)^{12} \left(\omega _1^2-4 \omega _3^2\right)}+\frac{135 q Q \sqrt{(1-\mu )^2} \sqrt{\mu ^2} \omega _1 \omega _3^2}{(1-\mu )^6 \mu ^6 \left(\omega
_1^2-4 \omega _3^2\right)}-}\\
 {\frac{270 q Q \sqrt{(1-\mu )^2} \sqrt{\mu ^2} \omega _1 \omega _3^2}{(1-\mu )^6 \mu ^5 \left(\omega _1^2-4 \omega _3^2\right)}+\frac{135 q Q
\sqrt{(1-\mu )^2} \sqrt{\mu ^2} \omega _1 \omega _3^2}{(1-\mu )^8 \mu ^4 \left(\omega _1^2-4 \omega _3^2\right)}+}\\
 {\frac{135 q Q \sqrt{(1-\mu )^2} \sqrt{\mu ^2} \omega _1 \omega _3^2}{2 (1-\mu )^7 \mu ^4 \left(\omega _1^2-4 \omega _3^2\right)}+\frac{135 q
Q \sqrt{(1-\mu )^2} \sqrt{\mu ^2} \omega _1 \omega _3^2}{(1-\mu )^6 \mu ^4 \left(\omega _1^2-4 \omega _3^2\right)}-}\\
 {\frac{270 q Q \sqrt{(1-\mu )^2} \sqrt{\mu ^2} \omega _1 \omega _3^2}{(1-\mu )^8 \mu ^3 \left(\omega _1^2-4 \omega _3^2\right)}-\frac{135 q Q
\sqrt{(1-\mu )^2} \sqrt{\mu ^2} \omega _1 \omega _3^2}{2 (1-\mu )^7 \mu ^3 \left(\omega _1^2-4 \omega _3^2\right)}+}\\
 {\frac{135 q Q \sqrt{(1-\mu )^2} \sqrt{\mu ^2} \omega _1 \omega _3^2}{(1-\mu )^8 \mu ^2 \left(\omega _1^2-4 \omega _3^2\right)}+\frac{270 q^2
\omega _1^2 \omega _3}{(1-\mu )^{10} \left(-4 \omega _1^2+\omega _3^2\right)}+\frac{270 Q^2 \omega _1^2 \omega _3}{\mu ^8 \left(-4 \omega _1^2+\omega
_3^2\right)}-}\\
 {\frac{540 q^2 \mu  \omega _1^2 \omega _3}{(1-\mu )^{10} \left(-4 \omega _1^2+\omega _3^2\right)}+\frac{270 q^2 \mu ^2 \omega _1^2 \omega _3}{(1-\mu
)^{10} \left(-4 \omega _1^2+\omega _3^2\right)}+\frac{540 q Q \sqrt{(1-\mu )^2} \sqrt{\mu ^2} \omega _1^2 \omega _3}{(1-\mu )^6 \mu ^5 \left(-4 \omega
_1^2+\omega _3^2\right)}-}\\
 {\frac{540 q Q \sqrt{(1-\mu )^2} \sqrt{\mu ^2} \omega _1^2 \omega _3}{(1-\mu )^6 \mu ^4 \left(-4 \omega _1^2+\omega _3^2\right)}+\frac{54 q^2
\omega _1 \omega _3^2}{(1-\mu )^{10} \left(-4 \omega _1^2+\omega _3^2\right)}+\frac{54 Q^2 \omega _1 \omega _3^2}{\mu ^8 \left(-4 \omega _1^2+\omega
_3^2\right)}-}\\
 {\frac{108 q^2 \mu  \omega _1 \omega _3^2}{(1-\mu )^{10} \left(-4 \omega _1^2+\omega _3^2\right)}+\frac{54 q^2 \mu ^2 \omega _1 \omega _3^2}{(1-\mu
)^{10} \left(-4 \omega _1^2+\omega _3^2\right)}+\frac{108 q Q \sqrt{(1-\mu )^2} \sqrt{\mu ^2} \omega _1 \omega _3^2}{(1-\mu )^6 \mu ^5 \left(-4 \omega
_1^2+\omega _3^2\right)}-}\\
 {\frac{108 q Q \sqrt{(1-\mu )^2} \sqrt{\mu ^2} \omega _1 \omega _3^2}{(1-\mu )^6 \mu ^4 \left(-4 \omega _1^2+\omega _3^2\right)}-\frac{54 q^2
\omega _3^3}{(1-\mu )^{10} \left(-4 \omega _1^2+\omega _3^2\right)}-}\\
 {\frac{54 Q^2 \omega _3^3}{\mu ^8 \left(-4 \omega _1^2+\omega _3^2\right)}+\frac{108 q^2 \mu  \omega _3^3}{(1-\mu )^{10} \left(-4 \omega _1^2+\omega
_3^2\right)}-\frac{54 q^2 \mu ^2 \omega _3^3}{(1-\mu )^{10} \left(-4 \omega _1^2+\omega _3^2\right)}-}\\
 {\left.\frac{108 q Q \sqrt{(1-\mu )^2} \sqrt{\mu ^2} \omega _3^3}{(1-\mu )^6 \mu ^5 \left(-4 \omega _1^2+\omega _3^2\right)}+\frac{108 q Q \sqrt{(1-\mu
)^2} \sqrt{\mu ^2} \omega _3^3}{(1-\mu )^6 \mu ^4 \left(-4 \omega _1^2+\omega _3^2\right)}\right);}\)

\noindent\( {\text{when}\text{  }\mu \to 0.00025,A\to 0.00025,q\to 0.025\text{  }\text{and} Q\text{-$>$}0.00025}\)

\noindent\( {D_2=-3.59996\times 10^{15} \left(-30720.2+\omega _1\right) \omega _1-\frac{8.64\times 10^{18} \omega _1^2}{\omega _3}+}\\
 {\left(-1.12895\times 10^{20}+1.72798\times 10^{16} \omega _1\right) \omega _3+\frac{\left(2.55997\times 10^{19}-5.7599\times 10^{15} \omega
_1\right) \omega _3^2}{\omega _1}+}\\
 {\frac{\omega _1^2 \left(2.21184\times 10^{20} \omega _1+5.5296\times 10^{19} \omega _3\right)}{-4. \omega _1^2+1. \omega _3^2}+\frac{\omega
_1^2 \left(-4.31996\times 10^{18} \omega _1-4.31996\times 10^{18} \omega _3\right)}{-0.25 \omega _1^2+1. \omega _3^2}}\)

\noindent\( {\text{when}\text{  }\mu \to 0.00025,A\to 0.00025,q\to 0.025,Q\text{-$>$}0.00025 \text{and} \omega _3\to 1 }\)
Using   all the values of coefficients  and algebraic manipulations  we get the following expression for $D_2$:\\
\noindent\( {D_2=5.76096\times 10^{14}-\frac{3.2\times 10^{14}}{\omega _1}-1.92\times 10^{11} \omega _1+}\\
 {3.6\times 10^{10} \omega _1^2+\frac{1.152\times 10^{15} \omega _1}{-4+\omega _1^2}-\frac{1.44\times 10^{14} \omega _1^2}{-4+\omega _1^2}-\frac{4.32\times
10^{14} \omega _1^3}{-4+\omega _1^2}+}\\
 {0.00025 \left(-2.30423\times 10^{23}+\frac{1.024\times 10^{23}}{\omega _1}+8.30131\times 10^{22} \omega _1-\right.}\\
 {3.45744\times 10^{22} \omega _1^2-\frac{2.21184\times 10^{23}}{1-4 \omega _1^2}+\frac{2.21184\times 10^{23} \omega _1}{1-4 \omega _1^2}+}\\
 {\left.\frac{1.10592\times 10^{24} \omega _1^2}{1-4 \omega _1^2}-\frac{5.5296\times 10^{23} \omega _1}{-4+\omega _1^2}+\frac{6.912\times 10^{22}
\omega _1^2}{-4+\omega _1^2}+\frac{2.0736\times 10^{23} \omega _1^3}{-4+\omega _1^2}\right)}\)

\begin{figure}[h]
 \includegraphics{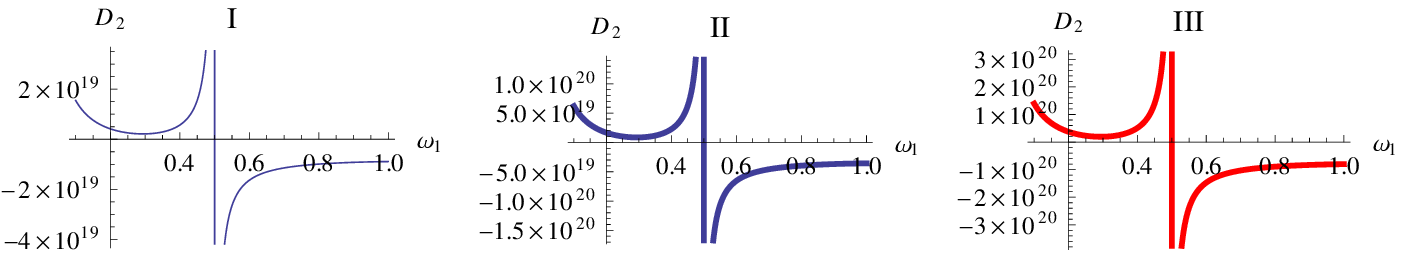}
\includegraphics{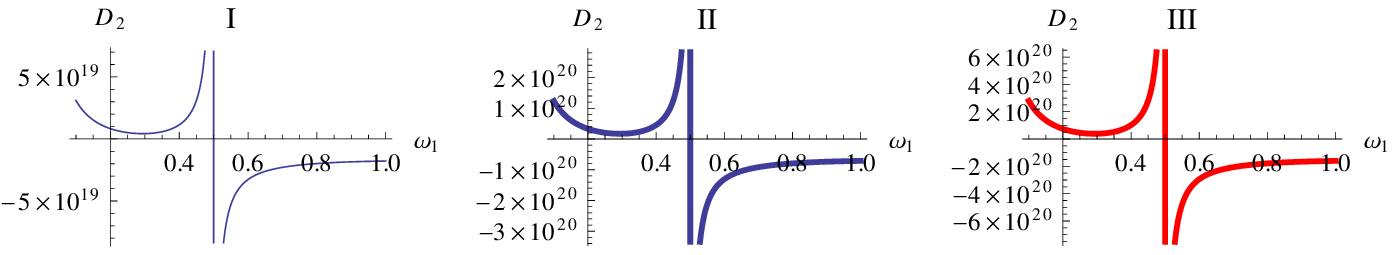}
\caption{$D_2$ $\text{when}\text{  }\mu \to 0.00025,A\to 0.00025,q_2 \to 0.25,0.50,$ ~\text{and}~$ 0.75, Q\to 0.00025 \text{ }  \text{and} \text{ } \omega _3\to 1 $\label{fig:fig1}}
\end{figure}

\begin{figure}[h]
\includegraphics{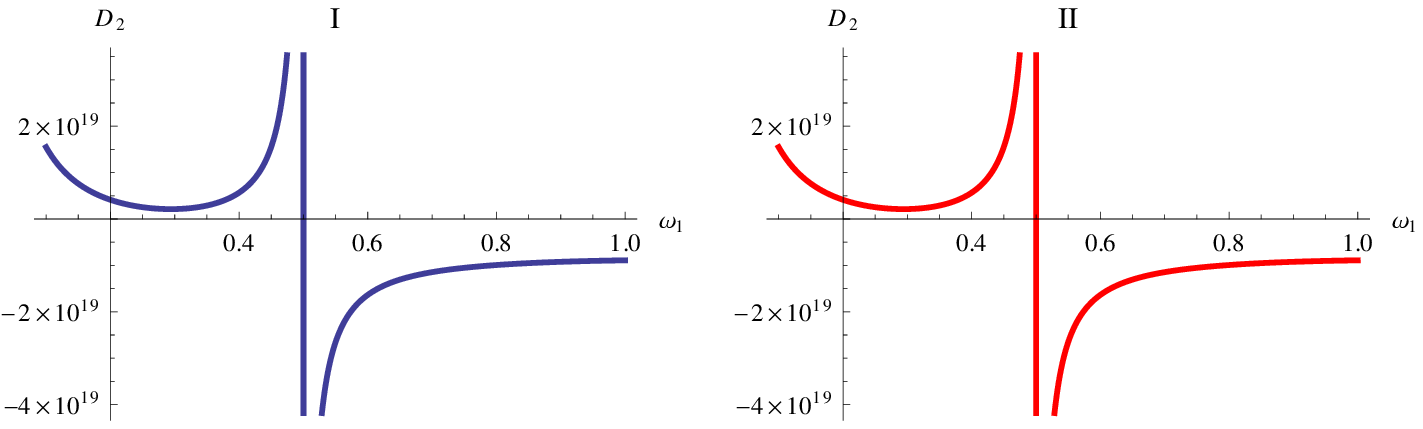}
\includegraphics{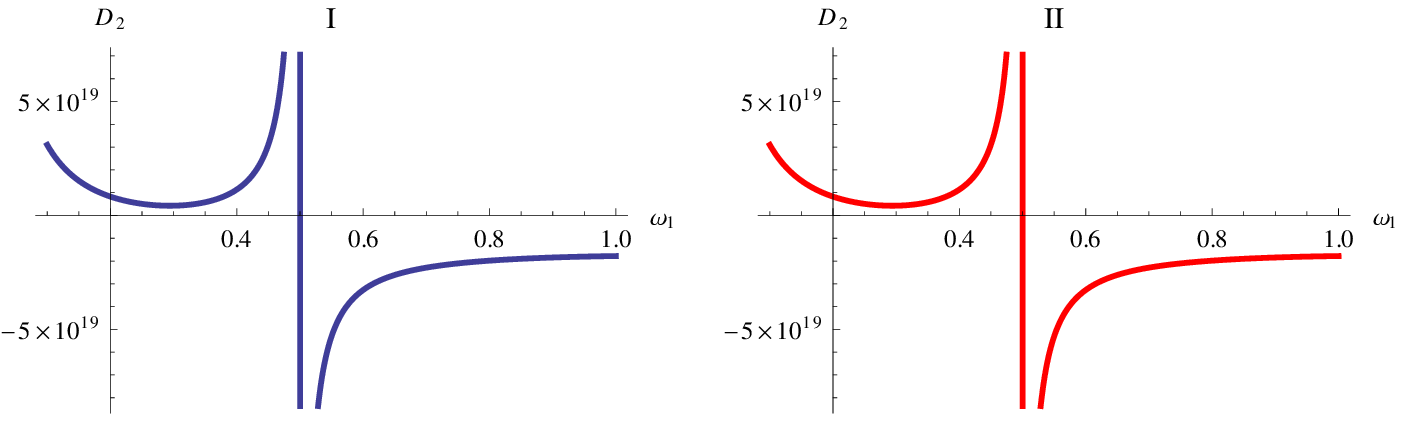}
\caption{\( {\text{when}\text{  }\mu \to 0.00025,A\to 0.00025,q_1 \to 0.50 ~ \text{and}~ 0.75,Q \to 0.00025 ~ \text{and}~ \text{  } \omega _3\to 1 }\)\label{fig:fig2}}
\end{figure}

 We have computed the value of $D_2$ numerically for various values of parameters when $\mu=0.0025$ and $\omega_3=1$.The graphs are plotted $D_2$ versus $\omega_1$.
In figures \ref{fig:fig1},   plot I,II,III represent the respective values of $q_2=0.25,0.50$, and $0.75$ respectively when $q_1=0.25$ fixed and the  vertical line shown in each graph is an asymptote.
Effect of $A_2$ in which first three curves are for $A_2=0.0025$ and second three curves belong to $A_2=0.0050$.
Similarly we have also obtained the effect of $q_1$ in figure \ref{fig:fig2} in which plot I and II correspond to $q_1=0.50$ and $0.75$ respectively when $q_2=0.25$ fixed. Here first  two curves are plotted for $A_2=.0025$ and second two curves for $A_2=.0050.$

\section{Conclusion}
It is evident from all the above figures that the curves are in rectangular hyperbolic forms with singularity at $\omega_1=.50$ and we find that $D_2\neq 0$  for various values of different parameters.
Thus we conclude that, according to Arnold's theorm since $D_2\neq 0$,out of plane equilibrium point $L_6$ is stable in non-linear sense.

{\bf Acknowledgements:} We are thankful to D.S.T. Govt. of India,New Delhi for sanctioning a project SR/S4/MS:380/06,13/5/2008.
We are also thankful to Dr Badam Singh  Kushvah Department of Applied Mathematics, ISM, Dhanbad, India for valuable suggestions in preparing this manuscript during our visit to IUCAA, Pune.
\bibliographystyle{spbasic}




\end{document}